\documentclass[12pt, reqno]{amsart}

\usepackage{graphicx,color,amsmath,amssymb,mathrsfs}
\usepackage{amsfonts,amsthm,epsfig,epstopdf,hyperref, mathtools}

\newtheorem{theorem}{Theorem}[section]

\newtheorem{corollary}[theorem]{Corollary}

\newcommand{\R}{{\mathbb{R}}}

\newcommand{\M}{{\mathcal{M}}}

\newcommand{\SSS}{{S_{p,n}}}

\newcommand{\vertiii}[1]{{\left\vert\kern-0.25ex\left\vert\kern-0.25ex\left\vert #1 
    \right\vert\kern-0.25ex\right\vert\kern-0.25ex\right\vert}}

\newcommand{\na}{{\nabla}}

\numberwithin{equation}{section}
\begin{document}

\title[A note on strong-form stability for the Sobolev inequality]{A note on strong-form stability for the Sobolev inequality}


\author{Robin Neumayer}
\address{School of Mathematics,
				Institute for Advanced Study, Princeton, NJ}
\email{neumayer@ias.edu}

\begin{abstract}
In this note, we establish a strong form of the quantitive Sobolev inequality in Euclidean space for $p \in (1,n)$. Given any function $u \in \dot W^{1,p}(\mathbb{R}^n)$, the gap in the Sobolev inequality controls $\| \nabla u -\nabla v\|_{p}$, where $v$ is an extremal function for the Sobolev inequality.
  \end{abstract}

\maketitle

\section{Introduction}
Sobolev inequalities, broadly speaking, establish integrability or regularity properties of a function in terms of the integrability of its gradient. A fundamental example is the classical Sobolev inequality on Euclidean space, which states the following. Given $n\geq 2$ and $p \in (1,n)$, there exists a constant $S=S(n,p)$ such that 
 \begin{equation}\label{Sobolev} \| \na u \|_{p} \geq S \|u\|_{{p^*}}. \end{equation}
for any function $u \in \dot W^{1,p}(\R^n)$. Here, $p^* = np/(n-p)$, and $\dot W^{1,p}(\R^n)$ is the space of functions such that $u \in L^{p^*}(\R^n)$ and $|\na u|\in L^p(\R^n)$.
Let us take $S$ to be the largest possible constant for which \eqref{Sobolev} holds. Aubin \cite{aubin1976} and Talenti \cite{talenti1976best} determined that equality is achieved in \eqref{Sobolev} for the function 
\begin{equation*}\label{v1}
\bar v(x) = {\left(1+|x|^{p'}\right)^{(p-n)/p}},
\end{equation*}
as well as its translations, dilations, and constant multiples. Here and in the sequel, we let $p' = p/(p-1)$ denote the H\"{o}lder conjugate of $p$. In fact, these functions are the only such extremal functions for \eqref{Sobolev}, and we will let
\[
\M = \Big\{ v\  \big| \ v(x) = c\, \bar v\left(\lambda(x-y)\right) \text{ for some } c \in \R,\,\lambda\in\R_+,\,y\in\R^n \Big\}
 \]
denote this $(n+2)$-dimensional space of extremal functions. 

Brezis and Lieb raised the question of quantitative stability for the Sobolev inequality in \cite{brezis1985sobolev}, asking whether the deviation of a given function from attaining equality in \eqref{Sobolev} controls its distance to the family of extremal functions $\M$. The strongest notion of distance that one expects to control is the $L^p$ norm between gradients. With this in mind, let us define the {\it asymmetry} of a function $u \in \dot W^{1,p}(\R^n)$ by
\begin{align*}
A(u)  & = \inf\left\{ \frac{ \| \na u - \na v\|_{p}}{\|u\|_{p^*}}	 : v\in \M 
\right\}\,.
\end{align*}
Note that $A(u)$ is invariant under the symmetries of the Sobolev inequality (translations, dilations, and constant multiples) and is equal to zero if and only if $u \in \M$.
	To quantify the deviation from equality in \eqref{Sobolev}, we define the {\it deficit} of a function $u \in \dot W^{1,p}(\R^n)$ to be 
\begin{align*}
\delta(u) &=	\frac{ \| \na u\|_{p}^{p'} - S^{p'}\|u\|_{{p^*}}^{p'} }{\| u\|_{p^*}^{p'}} \qquad  \qquad \text{ if } p < 2\\
\shortintertext{and}
 \delta(u) &= \frac{\| \na u\|_{p}^p - S^p\|u\|_{{p^*}}^p}{\| u\|_{p^*}^p} \qquad \qquad \ \ \ \text{ if } p \geq 2\,.
\end{align*}
Like the asymmetry, the deficit is a non-negative functional that is invariant under translations, dilations, and constant multiples, and is equal to zero if and only if $u\in \M$.

By way of a concentration compactness argument as in \cite{lions1985}, one readily establishes the qualitative stability of \eqref{Sobolev}. That is, if $\{u_i\}$ is a sequence of functions with $\delta(u_i) \to 0$, then $A(u_i)\to 0$.
The first quantitative result was established in the case $p=2$ in \cite{BianchiEgnell91}, where Bianchi and Egnell showed that there is a dimensional constant $C$ such that
\begin{align*}
	A(u)^2\leq C\,\delta(u)  \,.
\end{align*}
This result, in addition to being optimal in the strength of the distance controlled, is sharp in the sense that the exponent $2$ cannot be replaced by a smaller one. The proof relies strongly on the fact that $W^{1,2}(\R^n)$ is a Hilbert space, and in the absence of this structure, the case when $p\neq 2$ has proven much more difficult to treat. Nevertheless, in \cite{ciafusmag07}, Cianchi, Fusco, Maggi, and Pratelli  established a quantitative stability result in which the deficit controls the distance of a function to $\M$ in terms of the $L^{p^*}$ norm; see Theorem~\ref{CFMP Theorem} below for a precise statement. The argument combines symmetrization arguments in the spirit of \cite{FMP08} with  a mass transportation argument in one dimension. More recently, in \cite{FN19}, Figalli and the author strengthened this result in the case $p\geq 2$ by showing that the deficit of a function controls a power of $A(u)$. The main idea there was to view $W^{1,p}(\R^n)$ as a weighted Hilbert space and to establish a spectral gap for the linearized operator in the second variation as in \cite{BianchiEgnell91}. However, bounding the difference between the deficit and the second variation required the use of the main result of \cite{ciafusmag07}.

In this note, we establish a reduction theoreom that, paired with \cite{ciafusmag07}, allows us to deduce a strong-form quantitative stability result in which the deficit of a function controls a power of $A(u)$. For $p\geq 2$, this recovers the main result of \cite{FN19} with a simpler proof, while in the case $p \in (1,2)$, it provides the first known quantitative estimate for \eqref{Sobolev} at the level of gradients.

\begin{theorem}\label{MainThm} Fix $n\geq 2$ and $p \in (1,n).$
There exist constants $C_1(n,p)$ and $C_2(n,p)$ such that the following holds. For any $u \in \dot W^{1,p}(\R^n)$ and for any $v\in \M$ with $\| u\|_{p^*} = \| v\|_{p^*},$ we have  
\begin{equation} \label{eqn: main reduction}
 \left(\frac{\| \na u - \na v\|_{p}}{\| u\|_{p^*}}\right)^{\alpha} \leq C_1\, \delta(u) +C_2 \frac{ \|u-v\|_{{p^*}}}{\| u\|_{p^*}}\,.
 \end{equation}
 Here, $\alpha = p'$ if $p\in (1,2)$ and $\alpha = p$ if $p \in [2,n)$.
\end{theorem}
Pairing Theorem~\ref{MainThm} with the main result of \cite{ciafusmag07} (Theorem~\ref{CFMP} below), we establish the following quantitative estimate.
\begin{corollary}\label{MainCor} Fix $n \geq 2$ and $p \in (1,n)$. There exist constants $C= C(n,p)$ and $\beta = \beta(n,p)$ such that the following holds. For any $u \in \dot W^{1,p}(\R^n)$, we have
\begin{equation} \label{fullstability2}
 A(u)^{\beta}  \leq C\,\delta(u)\,.
 \end{equation}
\end{corollary}
The value of $\beta$ in Corollary~\ref{MainCor} is given by
\[
\beta = \begin{cases}
	p'\left(p^*\left(3 + 4p -\frac{3p+1}{n}\right)\right)^2 & \text{ if }p \in (1,2)\\
	p\left(p^*\left(3 + 4p -\frac{3p+1}{n}\right)\right)^2 & \text{ if } p \in [2,n)\,.
\end{cases}
\]
 
The proof of Theorem~\ref{MainThm} is elementary and at its core relies on the convexity of the function $t\mapsto t^p$. It is inspired by the recent paper \cite{HyndSeuffert}, in which Hynd and Seuffert give a qualitative description of extremal functions in (a certain form of) Morrey's inequality. Interestingly, they are able to establish a quantitative stability result, even without knowing the explicit form of extremal functions. 

 Quantitive stability for Sobolev-type inequalities has been a topic of interest in recent years.  Closely related to the main results here, a strong-form quantitative stability result was shown for the Sobolev inequality \eqref{Sobolev} with $p=1$ in \cite{figmagpraa}, following \cite{FMPSobolev, cianchi06}. Quantitative stability results have also been shown for (a different form of) Morrey's inequality \cite{cianchi08}, the log-Sobolev inequality \cite{IndreiMarcon, Bobkovetal, fathi2014quantitative}, the higher order Sobolev inequality \cite{BaWeWi, GazzolaWeth}, the fractional Sobolev inequality \cite{ChenFrankWeth},  Gagliardo-Nirenberg-Sobolev inequalities \cite{CarlenFigalli,  DT2, DT, Seuffert2, NguyenGNS}, and Strichartz inequalities \cite{Negro}.
 
 More broadly, strong-form stability estimates (in which the gap in a given inequality controls the strongest possible norm, typically involving the oscillation of a set or function) have been studied for various functional and geometric inequalities. For instance, such results have been shown for isoperimetric inequalities in Euclidean space \cite{fuscojulin11}, on the sphere \cite{BDuzFus13}, and in hyperbolic space \cite{BDS2015}, as well as for anisotropic \cite{Neumayer16} and Gaussian \cite{Eldan2015, BarBraJul14} isoperimetric inequalities.

Apart from their innate interest from a variational perspective, quantitative stability estimates have found applications in the study of  geometric problems \cite{figallimaggi11, CicaleseSparado13, KM14} and PDE
 \cite{CarlenFigalli, DT}. Certain applications, such as those in \cite{FMM, CNT}, necessitate strong-form quantitative estimates of the type established here.\\

\noindent{\it Acknowledgments:} The author is supported by Grant No. DMS-1638352 at the Institute for Advanced Study.

\section{Proofs of Theorem~\ref{MainThm} and Corollary~\ref{MainCor}}
In the proof of Theorem~\ref{MainThm}, we will make use of the following version of Clarkson's inequalities for vector-valued functions, which state the following. Let $F,G:\R^n\to \R^n$ with $|F|,|G|\in L^p(\R^n)$. Then
\begin{align}\label{eqn: clarkson p<2}
\left\| \frac{F+G}{2}\right\|_p^{p'} + \left\| \frac{F-G}{2}\right\|_{p}^{p'} & \leq \left( \frac{1}{2} \| F\|_{p}^p + \frac{1}{2} \| G\|_p^p\right)^{p'/p}	
\shortintertext{if $p \in (1,2)$, and }
\label{eqn: clarkson p>2}
\left\| \frac{F+G}{2}\right\|_p^p + \left\| \frac{F-G}{2}\right\|_p^p & \leq \frac{1}{2} \| F\|_{p}^p + \frac{1}{2} \| G\|_{p}^p.	
\end{align}
if $p \geq 2$. These inequalities were shown for scalar- and complex-valued functions in \cite{Clarkson}, and were extended to functions mapping from $\R$ to $\R^n$ in \cite{Boas}. Though Clarkson's inequalities have been generalized in a number of directions, we could not locate a reference for the precise form of \eqref{eqn: clarkson p<2} and \eqref{eqn: clarkson p>2}, so in Section~\ref{sec: appendix} we prove \eqref{eqn: clarkson p>2} and show how to deduce \eqref{eqn: clarkson p<2} from its scalar-valued analogue.
\begin{proof}[Proof of Theorem~\ref{MainThm}]
	We first consider the case $p\in(1,2)$. Applying \eqref{eqn: clarkson p<2} with $F = \na u$ and $G = \na v$, we find that
	\begin{align}\label{eqn: d}
		\left\| \frac{\na u - \na v}{2} \right\|_{p}^{p'} &\leq \left( \frac{1}{2}\| \na u \|_p^p + \frac{1}{2} \| \na v\|_{p}^p \right)^{p'/p} - \left\| \frac{\na u + \na v}{2} \right\|_p^{p'}
      	\end{align}
Next, the Sobolev inequality \eqref{Sobolev} implies that 
\begin{align}\label{eqn: e}
\| \na v\|_{p}^p &\leq \| \na u\|_p^p,\\
\shortintertext{and}
\label{eqn: f} \left\| \na u + \na v\right\|_p^{p'} & \geq S^{p'} \left\| u+v\right\|_{p^*}^{p'}\,.
\end{align}
In \eqref{eqn: e} we have used the assumption that $\| u\|_{p*}= \|v\|_{p^*}$.
Together \eqref{eqn: d}, \eqref{eqn: e}, and \eqref{eqn: f} imply that 
\begin{align}\label{eqn: g}
\left\| \frac{\na u - \na v}{2} \right\|_{p}^{p'} & \leq \| \na u \|_p^{p'} -S^{p'}\left\|\frac{ u+v}{2}\right\|_{p^*}^{p'}\,. \\
\shortintertext{Finally, we claim that}
\label{eqn: c}
		\left\| \frac{u+v}{2}\right\|_{p^*}^{p'} &\geq \| u\|_{p^*}^{p'} - {p'}\|u\|_{p^*}^{{p'}-1} \left\| \frac{u-v}{2}\right\|_{p^*} \,.
\\
\shortintertext{Indeed, Minkowski's inequality implies that }
\label{eqn: a}
	\left\| \frac{u+v}{2} \right\|_{p^*}^{p'} & \geq \left( \|u\|_{p^*} - \left\| \frac{u-v}{2}\right\|_{p^*}\right)^{p'}\,.	
	\end{align}
	Then, convexity of the function $t\mapsto t^{p'}$ implies that 
	\begin{align}\label{eqn: b}
		 \left( \|u\|_{p^*} - \left\| \frac{u-v}{2}\right\|_{p^*}\right)^{p'} & \geq \| u\|_{p^*}^{{p'}} - {p'}\| u\|_{p^*}^{{p'}-1} \left\| \frac{u-v}{2}\right\|_{p^*}\,.
	\end{align}
Together \eqref{eqn: a} and \eqref{eqn: b}
 imply \eqref{eqn: c}.
Finally, combining \eqref{eqn: g} and \eqref{eqn: c} and dividing through by $\|u\|_{p^*}^{p'}$ establishes the proof of \eqref{eqn: main reduction} with $C_1 = 2^{p'}$ and $C_2 = p' 2^{p'-1}$.

Next, the proof for the case $p\geq 2$ is completely analogous. 
Indeed, applying Clarkson's inequality \eqref{eqn: clarkson p>2} followed by the Sobolev inequality \eqref{Sobolev}, and then \eqref{eqn: c} (with $p$ replacing $p'$), we find that 
	\begin{align*}
 \left\|\frac{ \na u -\na v}{2}\right\|_{p}^p & \leq \frac{1}{2} \| \na u\|_{p}^p + \frac{1}{2} \| \na v\|_{p}^p -  \left\|\frac{ \na u + \na v}{2}\right\|_{p}^p\\
	& \leq \| \na u \|_p^p - S^p \left\|\frac{ u+v}{2}\right\|_{p^*}^p\\
	& \leq \| \na u \|_p^p - S^p\|u\|_{p^*}^p + p\|u\|_{p^*}^{p-1} \left\|\frac{ u-v}{2}\right\|_{p^*}\,.
	\end{align*}  
	Dividing by $\|u\|_{p^*}^p$ establishes \eqref{eqn: main reduction} with $C_1= 2^p$ and $C_2 = p2^{p-1}$.
\end{proof}
Now, let us recall the main result from \cite{ciafusmag07}. The notion of $L^{p^*}$ asymmetry considered there is
\[
\lambda(u) = \inf \left\{\frac{ \| u - v\|_{{p^*}}}{\|u\|_{{p^*}}} : v \in \M, \ \|v\|_{p^*} = \|u\|_{p^*} \right\} 
\]
\begin{theorem}[Cianchi, Fusco, Maggi, Pratelli] \label{CFMP Theorem} Fix $n \geq 2$ and $p \in (1,n)$. There exists a constant $C= C(n,p)$ such that the following holds. For any $u\in \dot{W}^{1,p}(\R^n)$,
\begin{equation}\label{CFMP} \lambda(u)^{\beta}\leq C\,\frac{ \|\na u \|_{p} - \SSS \| u \|_{p^*}}{\| u\|_{p^*}},\end{equation}
Here $\beta = \left(p^*\left(3 + 4p -\frac{3p+1}n\right)\right)^2$.
\end{theorem} 
We now prove Corollary~\ref{MainCor} by combining Theorems~\ref{MainThm} and \ref{CFMP Theorem}.
\begin{proof}[Proof of Corollary~\ref{MainCor}] 
The only point to check is that
\begin{equation}\label{eqn: order deficits}
\frac{ \| \na u\|_{p} -S\| u \|_{{p^*}}}{\|u\|_{{p^*} }}
 \leq \delta(u)\,. 
 \end{equation}
To see this, note that for any $q\geq 1$, the function $t\mapsto t^q-t$ is increasing for $t\geq 1.$ In particular, if $a \geq b \geq 1$, we have 
\begin{equation}\label{eqn: numbers}
a^q -b^q \geq  a-b\,.
\end{equation}
Let $a = \|\na u \|_{p} / S\|u \|_{p^*}$ and $b=1$. Then applying \eqref{eqn: numbers} with $q=p'$ for $p\in(1,2)$ and $q=p$ for $p\in [2,n)$ establishes \eqref{eqn: order deficits}. With this in hand, Corollary~\ref{MainCor} follows immediately from \eqref{eqn: main reduction} and \eqref{CFMP}.
\end{proof}

\section{Clarkson's inequalities for vector valued functions on $\R^n$}\label{sec: appendix}
For $p\in (1,2)$, Clarkson \cite{Clarkson} established the following inequality for, in particular, real numbers $a$ and $b$:
\begin{align}\label{eqn: clarkson for numbers}
	|a+b|^{p'} + |a-b|^{p'}  & \leq 2 (|a|^p + |b|^p)^{{p'}/p} \,.
\end{align}
 Let us see how to deduce \eqref{eqn: clarkson p<2} from \eqref{eqn: clarkson for numbers}.
We make use of the reverse Minkowski inequality: if $s\in(0,1)$, then for $(a_1,\dots, a_n)\subset \R^n$ and $(b_1,\dots, b_n)\subset \R^n$ we have 
\begin{align}\label{eqn: mink a}
\left( \sum |a_i|^s \right)^{1/s} + \left( \sum |b_i|^s \right)^{1/s} &\leq \left( \sum |a_i+b_i|^s \right)^{1/s}
\end{align}
This inequality follows from the concavity of the function $t\mapsto t^s.$ We take $s = 2/{p'}$ and let $a_i = |F_i+G_i|^{p'}$ and $b_i = |F_i-G_i|^{p'}$ for $i=1,\dots, n$. Here $F_i$ denotes the $i$th component of $F$ in some fixed basis. Then, applying \eqref{eqn: mink a} followed by \eqref{eqn: clarkson for numbers}, we find that
\begin{equation}\label{eqn: intermediate step}
\begin{split}
	|F-G|^{p'} + |F+G|^{p'} & \leq \left( \sum (|F_i+G_i|^{p'} + |F_i-G_i|^{p'})^{2/{p'}}\right)^{{p'}/2}	\\
&\leq 2\left( \sum (|F_i|^p + |G_i|^p)^{2/p}\right)^{{p'}/2}\,.
\end{split}	
\end{equation}
On the left-hand side, we have used $|F|$ to denote the Euclidean norm. Next, applying the usual form of Minkowski's inequality with $r = 2/p$ to $(|a_i|^p)$ and $(|b_i|^p)$, we find
\begin{align*}
	\left( \sum (|a_i|^p + |b_i|^p)^{2/p}\right)^{1/2} \leq \left( \left( \sum |a_i|^2\right)^{p/2} + \left(\sum|b_i|^2\right)^{p/2}\right)^{1/p}\,.
\end{align*}
 Pairing this with \eqref{eqn: intermediate step}, we find that 
\begin{equation}\label{eqn: intermediate step2}
|F+G|^{p'} + |F-G|^{p'}  \leq 2\left( 	|F|^p + |G|^p\right)^{{p'}/p}.
\end{equation}
Finally, we make use of the integral form of \eqref{eqn: mink a}: for $s \in (0,1)$, we have
\begin{align}\label{eqn: mink b}
	\| h_1\|_{L^s(\R^n)} + \|h_2\|_{L^s(\R^n)}& \leq \|h_1 + h_2\|_{L^s(\R^n)}\,.
\end{align}
We apply \eqref{eqn: mink b} with $s=p/{p'}$ and with $h_1 = |F+G|^{{p'}}$ and $h_2 = |F-G|^{{p'}}$, and then apply \eqref{eqn: intermediate step2}, in order to find that
\begin{align}
\|F+G\|_p^{p'} + \| F-G\|_p^{p'} &\leq \left( \int \left( |F+G|^{p'} + |F-G|^{p'} \right)^{p/{p'} } \right)^{{p'}/p}\\
& \leq 2 \left( \int |F|^p +|G|^p\right)^{{p'}/p}\,.
\end{align}
This establishes \eqref{eqn: clarkson p<2}. The corresponding inequality \eqref{eqn: clarkson p>2} for $p\geq 2$ is straightforward. 
Note that $a^{p/2} +b^{p/2} \leq (a+b)^{p/2}$ for $p\geq 2.$ Applying this property and then expanding the squares, we have 
\begin{equation}\label{eqn: 2a}
\begin{split}
	|F+G|^p + |F-G|^p & = \left( \sum |F_i+G_i|^2 \right)^{p/2} + \left( \sum |F_i-G_i|^2 \right)^{p/2}\\
	&\leq \left( \sum (|F_i +G_i|^2 + |F_i-G_i|^2 )\right)^{p/2}\\
	& = \left(2\left( |F|^2 + |G|^2 \right)\right)^{p/2}	\,.
\end{split}	
\end{equation}
Finally, convexity of the function $t\mapsto t^{p/2}$ implies that
\begin{equation}\label{eqn: 2b}
\begin{split}
	\left(2 \left(|F|^2 + |G|^2 \right)\right)^{p/2} & = 2^{p} \left(\frac{|F|^2}{2}+ \frac{|G|^2}{2} \right)^{p/2} \\
	& \leq 2^{p-1} (|F|^p + |G|^p).	
\end{split}	
\end{equation}
We combine \eqref{eqn: 2a} and \eqref{eqn: 2b} and integrate to conclude the proof of \eqref{eqn: clarkson p>2}.

 \bibliographystyle{alpha}
\bibliography{references4}

\end{document}